# A Challenging 7-Fold Tiling Puzzle

Theo P. Schaad

**Abstract:**

A quasiperiodic 7-fold rhombic tiling is constructed with an iterative substitution scheme. The inflation factor is δ=5.04892…, the square of the longer diagonal of a regular heptagon. There are many substitutions possible that fill larger similar tiles with three base shapes but finding the matching rules (how the larger tiles fit together to make even larger tiles) turns into a challenging puzzle. At the end, a new solution is found with seven substitution rules. The relationship to a higher-order Fibonacci series is explored.

Outline:

A set of tiles for covering a surface is composed of three types of tiles. The base shape of each one of them is a rhombus, with acute angles of π/7, 2π/7, and 3π/7 radians (A, B, and C tiles). Such a set is called 7-fold, as the B tiles can be arranged in 7-fold rotational symmetry forming a 7-pointed star. A substitution scheme is sought to fill similar larger tiles with the base tiles. The process can then be repeated, and an arbitrarily large surface can be tiled. Several 7-fold sets have been discovered with only one substitution rule for each tile. One particularly challenging puzzle has been to find a scheme where the edges of the larger tiles are inflated by δ=5.04892… Alexey Madison[14] found a solution with nine different substitutions. Here, the process of finding other solutions is described. Finally, a different pattern was discovered with seven substitutions.

Background:

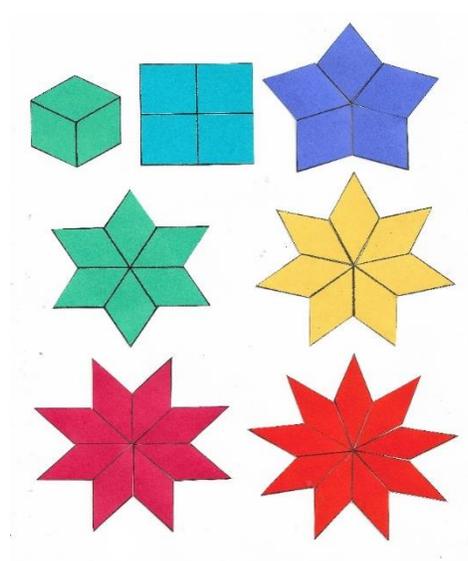

**Fig. 1:** N-fold stars play an important role in many n-fold tilings. The long diagonal of the 5-fold star rhombus is the golden ratio celebrated in the arts and mathematics. Likewise, the long diagonal of the 7-fold star rhombus has unique properties, somewhat overlooked.

Tile sets with n-fold symmetry have been studied extensively since the discovery of the Penrose 5-fold quasiperiodic pattern[1]. The 5-fold pattern is aesthetically pleasing, in part because it incorporates the golden ratio Φ. $Φ_n$ is the radius of a circle that surrounds the n-fold star (and is also a diagonal of the rhombus with an acute angle of $2π/n$).

$Φ_n = 2 \cos(π/n)$

| n-fold | $Φ_n$ | Value | a.k.a. | Ref. |
|---|---|---|---|---|
| 3-fold | 1 | 1 | | |
| 4-fold | √2 | 1.414 | Silver Ratio | |
| 5-fold | (1+√5)/2 | 1.618 | Golden Ratio | 1 |
| 6-fold | √3 | 1.732 | | 3 |
| 7-fold | 2 cos(π/7) | 1.802 | Magic 7-fold Φ | 2 |
| | | | Optimal Trisection | 6 |
| | | | Golden Trisection | 10 |
| 8-fold | 2 cos(π/8) | 1.848 | | 4,5 |
| 9-fold | 2 cos(π/9) | 1.879 | Quadrisection | 6 |
| 10-fold | 2 cos(π/10) | 1.902 | | 7 |
| 11-fold | 2 cos(π/11) | 1.919 | Pentasection | 6 |
| 12-fold | 2 cos(π/12) | 1.932 | | 8,9 |

In my previous study of 7-fold tilings[2], I emphasized the number $Φ_7$ (Φ from now on) which I called the 7-fold magic number. It is the longer diagonal of the 7-fold star rhombus (B tile) with edges of unit length. Following the work of P. Steinbach[6], I will now add a second important number to the 7-fold tilings. As shown in Figure 2, P. Steinbach noted that a regular heptagon (7-sided polygon) with sides of unit length has two diagonals, one shorter with length Φ, and another longer with length ρ. He discovered that the proportions 1: Φ: ρ have amazing properties and calls it the optimal trisection of a line segment. It is an extension of the golden bisection which leads to the golden ratio. It turns out that the trisection also shows up in the ratios of the areas of the three base tiles in 7-fold tilings.

Areas:

| Tile | Area | Value | Related to Φ, ρ |
|---|---|---|---|
| A | A = sin(π/7) | 0.434 | |
| B | B = sin(2π/7) | 0.782 | Φ*A |
| C | C = sin(3π/7) | 0.975 | ρ*A |

Area Ratios:

| Ratio | Value | Related to Φ, ρ |
|---|---|---|
| B/A | 1.802 | Φ |
| C/B | 1.247 | ρ-1 |
| C/A | 2.247 | ρ |

Diagonals:

| Tile | Diag. | Trig. | Value | Related to Φ, ρ |
|---|---|---|---|---|
| A | Short | $\gamma = 2 \sin(\pi/14)$ | 0.445 | ρ-Φ |
|   | Long | $\Gamma = 2 \cos(\pi/14)$ | 1.950 | $(\Phi+2)^{½}$ |
| B | Short | $\phi = 2 \sin(\pi/7)$ | 0.868 |  |
|   | Long | $\Phi = 2 \cos(\pi/7)$ | 1.802 | Φ |
| C | Short | $\psi = 2 \sin(3\pi/14)$ | 1.247 | ρ |
|   | Long | $\Psi = 2 \cos(3\pi/14)$ | 1.564 |  |

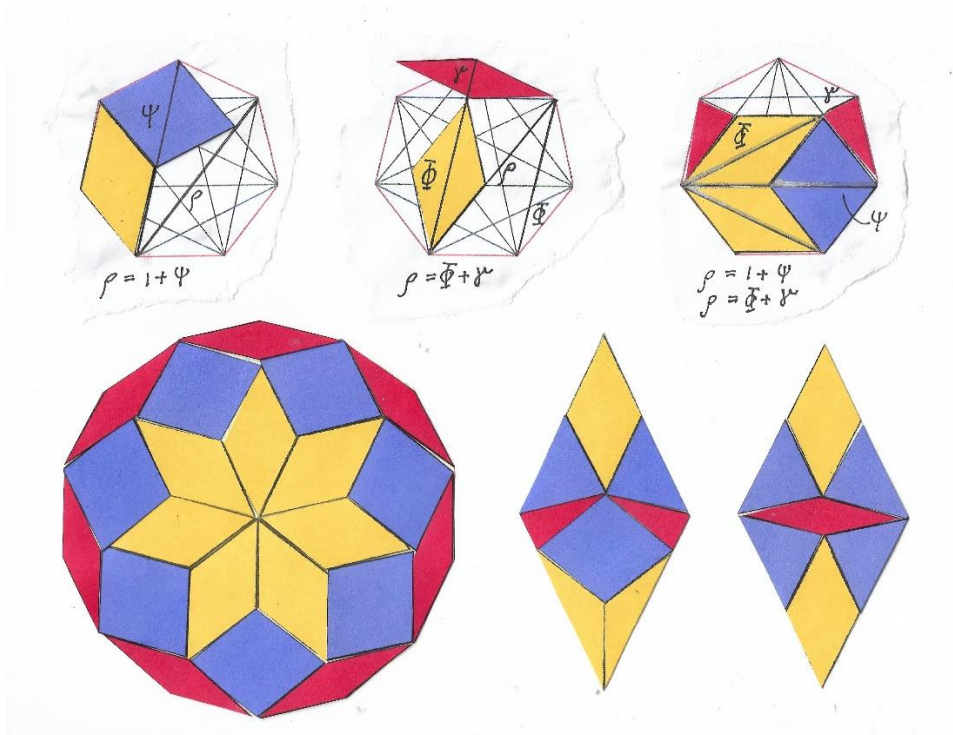

**Fig. 2:** A regular heptagon has a shorter diagonal Φ which is also the long diagonal of the yellow B tile. The longer diagonal ρ of the heptagon can be composed of a tile edge (unit length) plus the short diagonal ψ of the blue C tile, or alternatively, as the sum of Φ plus the short diagonal γ of the red A tile. The tetradecagon (14-sided regular polygon) is also an important motif of some 7-fold tilings. It can be filled with 7 each of the three base tiles. A circle drawn around it has a radius of ρ, the longer diagonal of the heptagon. A substitution tiling of the B tile is shown to the right with an inflation factor ρ (new edge length).

P. Steinbach found some remarkable properties in the trisection that he called "To Add is to Multiply". As shown in Fig. 2, it is possible to find a substitution rule for the B tile if the inflation factor is ρ (the longer diagonal of the heptagon or the radius around the tetradecagon in Fig. 2). It is the sum of an edge (unit length) and the short diagonal ψ of the blue C tile:

P = 1 + ψ

Now, if we were to inflate it again, we would multiply the inflation factor as ρ*ρ and end up with more tiles along the new inflated edge, which is an addition:

ρ*ρ = 1 + Φ + ρ

Although P. Steinbach stated his result in artistic terms, a Greek idea that allows similar figures to be arranged into another set of different proportions, it is the perfect description of a 7-fold quasiperiodic tiling. On the left is the inflation factor ρ*ρ = 5.04892… and on the right are the matching rules of the edges, which are a combination of rhomb sides (unit length) and three diagonals: the short diagonals of the skinny A rhomb and the thick C rhomb, and the long diagonal of the B tile.

I can now restate the goal of the puzzle: What are the substitution and matching rules for a 7-fold tiling with an inflation factor of δ = ρ$^2$ (= 5.04892…)? Previously[2], I had grouped all those solutions as Pattern #5 (an arbitrary index). I had found substitutions, but no matching rules for a set of only three base tiles.

I also showed that the substitution rule can be given as a matrix:

$A_1$          $A_0$

$B_1$ =   M    $B_0$

$C_1$          $C_0$

where

$$M = \begin{vmatrix} a & b & c \\ b & (a+c) & (b+c) \\ c & (b+c) & (a+b+c) \end{vmatrix}$$

Pattern #5 corresponds to the group (a,b,c)=(3,5,6), hence

$$M = \begin{vmatrix} 3 & 5 & 6 \\ 5 & 9 & 11 \\ 6 & 11 & 14 \end{vmatrix}$$

In words, a larger similar skinny $A_1$ rhomb is filled with 3 skinny $A_0$ rhombs, 5 $B_0$ rhombs, and 6 thick $C_0$ rhombs. So, we know exactly how many base tiles are needed for the inflated tiles. We just do not know the matching rules.

Penrose and Fibonacci Revisited:

To better understand Steinbach's optimal trisection and how it relates to the 7-fold pattern, it helps to look at the golden bisection and how it relates to the pentagonal Penrose pattern.

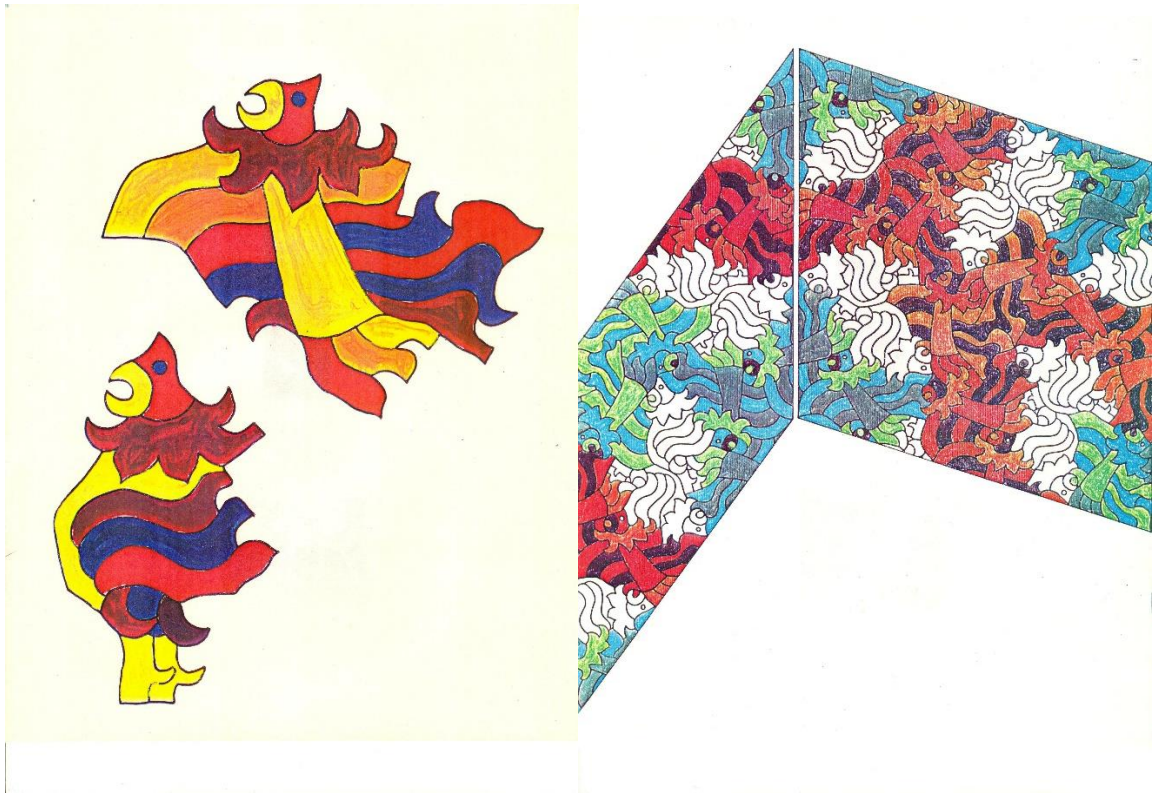

**Fig. 3:** I tiled Escheresque Penrose birds (left) into inflated similar rhombs (right). Here, as elsewhere, I cut the rhombs (birds) into half-rhombs (partial birds), with the caveat that the matching rules guarantee that the rhombs (birds) are fully reconstructed in the final tiling. Alternatively, one could use entire birds wherever the head is seen and eliminate those without a head. How many birds are in the inflated rhombs? Hint: Count heads! In the smaller tile, there are 5 small and 8 large birds, in the larger tile, 8 small and 13 large birds. Checking with the Penrose matrices in the text, it is the third generation (n=3). The bird's sizes (areas) are in the golden ratio. This decoration is from 1979.

A bisection of a line into two segments a, b (here a < b) with the golden ratio b: a = (a+b): b can be generated with

b       =       M       a

a+b                     b

where the matrix M is

M = $\begin{vmatrix} 0 & 1 \\ 1 & 1 \end{vmatrix}$

Powers lead to a series $M^n$

$\begin{vmatrix} 0 & 1 \\ 1 & 1 \end{vmatrix}$ $\begin{vmatrix} 1 & 1 \\ 1 & 2 \end{vmatrix}$ $\begin{vmatrix} 1 & 2 \\ 2 & 3 \end{vmatrix}$ $\begin{vmatrix} 2 & 3 \\ 3 & 5 \end{vmatrix}$ $\begin{vmatrix} 3 & 5 \\ 5 & 8 \end{vmatrix}$ $\begin{vmatrix} 5 & 8 \\ 8 & 13 \end{vmatrix}$ etc.

This is the Fibonacci series with

$\begin{vmatrix} 0 & 1 \\ 1 & 1 \end{vmatrix}^n = \begin{vmatrix} F_{n-1} & F_n \\ F_n & F_{n+1} \end{vmatrix}$

where $F_0 = 0$, $F_1 = 1$, and $F_n = F_{n-1} + F_{n-2}$ ($n \geq 2$)

$F_{n=0,1,2...} = 0, 1, 1, 2, 3, 5, 8, 13, ...$

For large n, the powers of the matrix converge to

$\begin{vmatrix} a: & b \\ b: & (a+b) \end{vmatrix} = \begin{vmatrix} 1: & \Phi \\ \Phi: & (\Phi+1) \end{vmatrix}$ (normalized to 1)

i.e., adjacent Fibonacci numbers approximate the golden ratio $\Phi$ for large n.

The substitution matrix for the Penrose pattern is

M = $\begin{vmatrix} 1 & 1 \\ 1 & 2 \end{vmatrix}$

$M^n = \begin{vmatrix} F_{2n-1} & F_{2n} \\ F_{2n} & F_{2n+1} \end{vmatrix}$

i.e., the Penrose matrices are the n even matrices of the Fibonacci matrices. For large n, both series converge to b: a = $\Phi$, the golden ratio, either as a ratio of line segments, or, in the Penrose tiling, the ratio of the areas of the two base tiles.

## The Third-Order Fibonacci Series and 7-Fold Tilings:

B. Franco[12,13] made the connection of a higher-order Fibonacci series to a heptagonal tiling that involves 5 base tiles, triangles, not rhombs. The third-order Fibonacci sequence can be defined as

$S_n = 2S_{n-1} + S_{n-2} - S_{n-3}$

This definition is linked to a paper of H. Terauchi[13] et al. The sequence is

$S_{n=0,1,2...}$ = 0, 0, 1, 2, 5, 11, 25, 56, 126, 636, 1429, …

It describes a transformation

$$\begin{pmatrix} C \\ B+C \\ A+B+C \end{pmatrix} = M \begin{pmatrix} A \\ B \\ C \end{pmatrix}$$

where

$$M = \begin{pmatrix} 0 & 0 & 1 \\ 0 & 1 & 1 \\ 1 & 1 & 1 \end{pmatrix}$$

Note how this matrix is a simple extension to a 3x3 matrix from the 2x2 Fibonacci matrix.

The powers are

$$M^n = \begin{pmatrix} S_n - S_{n-1} & S_n & S_{n+1} - S_n \\ S_n & S_{n+1} + S_{n-1} & S_{n+1} \\ S_{n+1} - S_n & S_{n+1} & S_n + S_{n+1} - S_{n-1} \end{pmatrix}$$

P. Steinbach calculated these matrices (n=1,2,3,4,5)

$$\begin{pmatrix} 0 & 0 & 1 \\ 0 & 1 & 1 \\ 1 & 1 & 1 \end{pmatrix} \quad \begin{pmatrix} 1 & 1 & 1 \\ 1 & 2 & 2 \\ 1 & 2 & 3 \end{pmatrix} \quad \begin{pmatrix} 1 & 2 & 3 \\ 2 & 4 & 5 \\ 3 & 5 & 6 \end{pmatrix} \quad \begin{pmatrix} 3 & 5 & 6 \\ 5 & 9 & 11 \\ 6 & 11 & 14 \end{pmatrix} \quad \begin{pmatrix} 6 & 11 & 14 \\ 11 & 20 & 25 \\ 14 & 25 & 31 \end{pmatrix}$$

Et voila! We immediately recognize that the 4th matrix is the substitution of the 7-fold Pattern #5, the object of this paper!

The series converges to Steinbach's optimal trisection ratios (normalized to 1)

1:     Φ:     ρ

Φ:     1+ρ:   Φ+ρ

ρ:     Φ+ρ:   1+Φ+ρ

Checking the 10[th] power, we get Φ ≈ 636/353 (≈ 1.802) and ρ ≈ 793/353 (≈ 2.246)

A Detour: A 7-Fold Tiling with Inflation Factor ρ

Analogous to the Penrose pattern that starts with the second Fibonacci matrix, it would be interesting to explore the second matrix of the third-order Fibonacci sequence. A start was already given in Figure 2 for the B tile. It is indeed possible to find substitutions for the A and C tiles as shown in Figure 4. But there are problems for the B tile in the second generation and the scheme does not work.

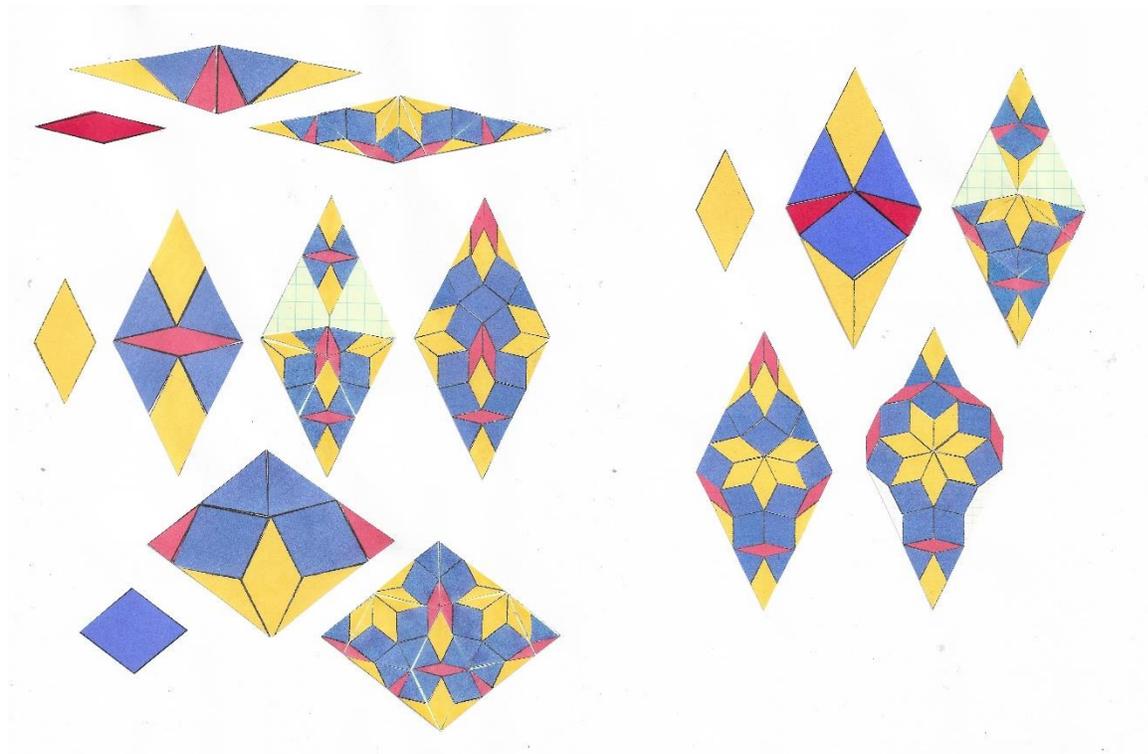

**Fig. 4:** A substitution scheme for the three base tiles with an inflation factor δ = ρ (the longer diagonal of a regular heptagon). There are solutions for the first generation, but problems arise trying to find the second generation. A rearrangement of a few base tiles leads to inflated tiles with an inflation factor δ = ρ². The last substitution features an example where it was not possible to use half-rhombs at the edges. The substitution is attractive with a 7-fold star, but it could never be assembled into a similar larger tile in the next generation as the tile is not 7-fold rotationally symmetric at the acute vertex.

## Substitution Rules for Pattern #5 with Inflation Factor $\rho^2$:

As shown in the introduction, P. Steinbach expressed the inflation factor $\delta = \rho^2$ (= 5.04892…) as

$$\rho^2 = 1 + \Phi + \rho$$

where $\rho = 1 + \psi$ or $\rho = \Phi + \gamma$. $\gamma$ is the short diagonal of the A tile, $\Phi$ the long diagonal of the B tile, and $\psi$ the short diagonal of the C tile. In the context of tiling, the inflation factor is the new side of the inflated tile. Limiting our choices only to inflations that tile edge-to-edge and use half-rhombs, the inflation factor could be any commutative arrangement of the addition

$$\delta = 1 + 1 + \Phi + \psi$$

$$\delta = 1 + \Phi + \Phi + \gamma$$

There is no symmetric solution, and the inflated edges are directional. Adjacent tiles can only match if the sequence is the same in the same direction.

Figures 6-8 show many examples of substitutions that are possible. Each tile shape contains the same number of base tiles as it is determined by the substitution matrix

$$M = \begin{vmatrix} 3 & 5 & 6 \\ 5 & 9 & 11 \\ 6 & 11 & 14 \end{vmatrix}$$

It is easiest to start finding substitutions for the skinny A rhomb. Preferable substitutions should be bilateral along the short diagonal. There are six possibilities. The inflation is (starting at the acute angle)

$$\delta = 1 + \Phi + \psi + 1 \text{ (A}_1 \text{ tile, all edges)}$$

$$\delta = \Phi + 1 + \gamma + \Phi \text{ (A}_2 \text{ tile, upper edges)}$$

$$\delta = 1 + \Phi + \psi + 1 \text{ (A}_2 \text{ tile, lower edges)}$$

In all cases, the areas near the vertices can only be tiled with an A or B tile, which then limits the choices for matching substitutions of the B and C tiles.

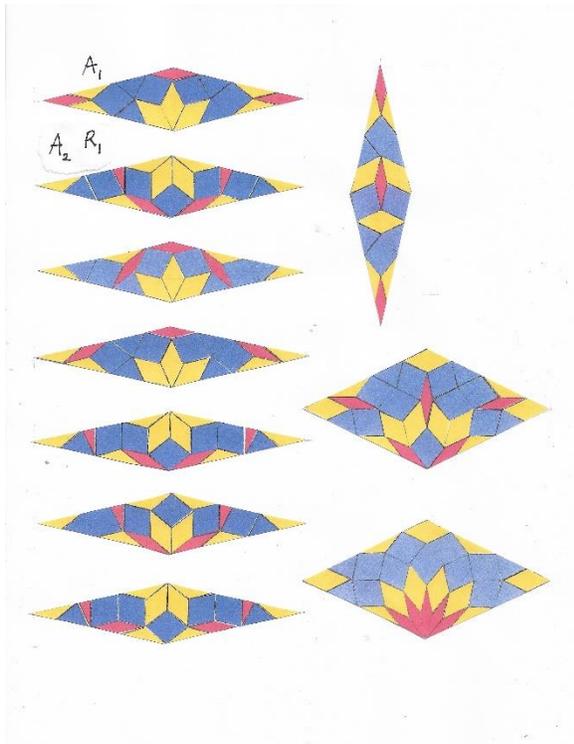

**Fig. 5:** Pattern #5 substitutions for the A file (skinny rhomb). Most are bilateral along the short diagonal; there is only one substitution that is bilateral along the long diagonal. $R_1$ is used in the Madison tiling, $A_1$ and $A_2$ in the Pattern #5 tiling. Two substitutions for the B tile along the short diagonal (lower right). The last one was nicknamed "sunrise", but, so far, has not led to a substitution tiling.

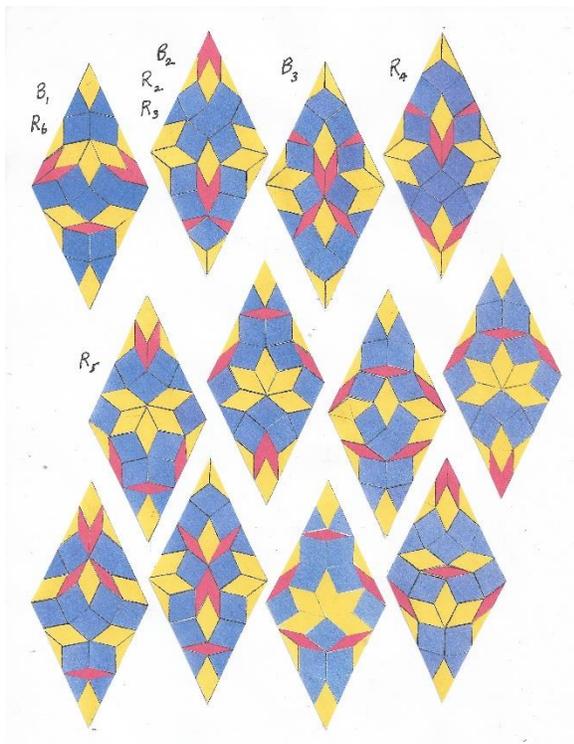

**Fig. 6:** Pattern #5 substitutions for the B tile. $R_2$ - $R_6$ are used in the Madison tiling, $B_1$ - $B_3$ in the Pattern #5 tiling. Only bilateral substitutions along the long diagonal are shown.

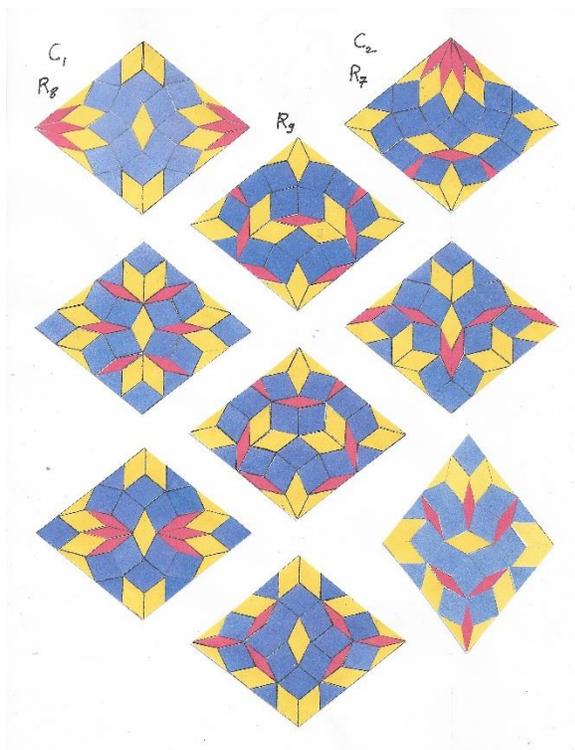

**Fig. 7:** Pattern #5 substitutions for the C tile. $R_7$ - $R_9$ are used in the Madison tiling, $C_1$ and $C_2$ in the Pattern #5 tiling. Four tiles have two mirror planes, but if they were inflated to the next generation, the central yellow B tile would break the symmetry, as it is directional (see examples in Fig. 6). Five tiles are bilateral along the short diagonal, only the last example is bilateral along the long diagonal only.

The inflated tiles are now all potential puzzle pieces that, hopefully, could be used to create larger similar tiles in the second generation. The process is near impossible to describe and is best performed by copying and cutting out the pieces. Furthermore, if more than one substitution is used for the same base shape, they need to be labeled or colored differently. I do not believe that a substitution scheme with only three tiles is possible, and my original search was in vain.

But there is at least one solution, and another one is shown later. Alexey Madison came up with an ingenious pattern that is arguably the most aesthetic 7-fold tiling discovered so far.

## Madison's 7-Fold Tiling with Nine Substitutions:

Alexey Madison's 7-fold tiling is described in a refereed paper in Structural Chemistry[14] (2018) that should be consulted for more information. It is also included in the Bielefeld Tiling Encyclopedia[15], but only with sparse information.

The substitution rules are shown in Figure 8 to give an overall impression of the substitutions without going into the details of the matching rules. Madison also recognized that the inflation factor is the square of the longer diagonal of a regular heptagon, based on Pautze[16].

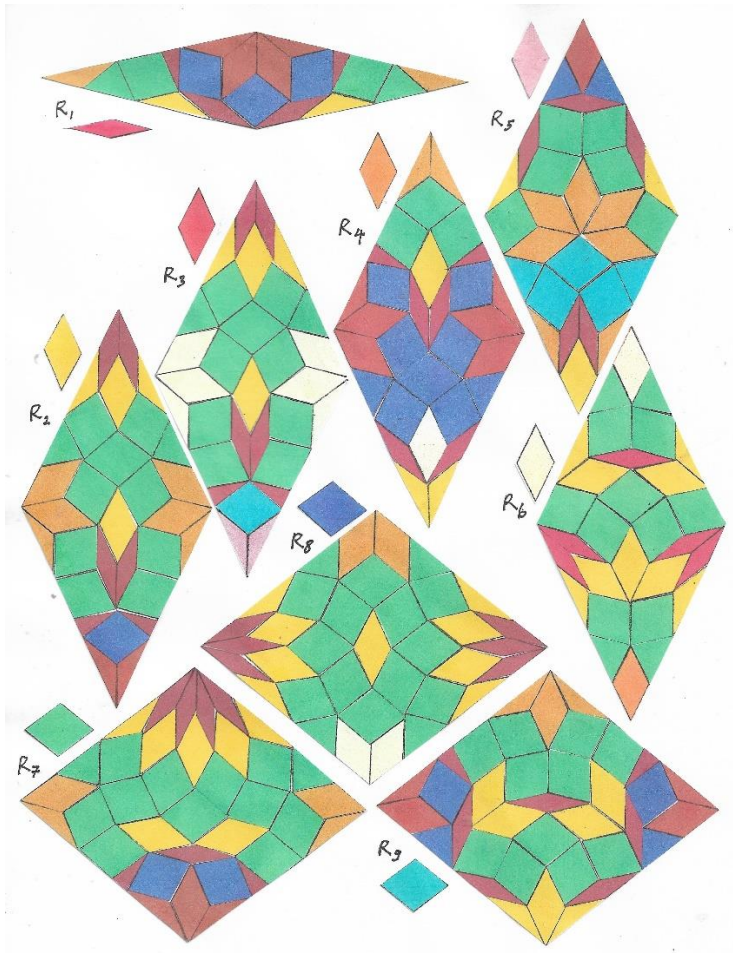

**Fig. 8:** Substitution rules for Madison's 7-fold tiling with nine substitutions. Colors are used to distinguish base tiles with different substitutions.

Nine rhombi are used to generate an infinite self-similar tiling. The most common tile is the green $R_7$ tile that creates a green background for 7- and 14-fold stars in a larger patch of the tiling (see Fig. 10). Most vertices of the inflated tiles display a portion of a star. Madison labels the vertices with $A_n$, $B_n$ (n=1,2,3), and C, not shown here. The matching rules require that an edge of adjacent tiles only fit with the same vertex labels, often completing a star at each vertex. In a larger patch, each of the five B tiles appear in a star configuration, thus the tiling has five 7-fold stars in different colors, which gives it its unique appeal as shown in Figure 9.

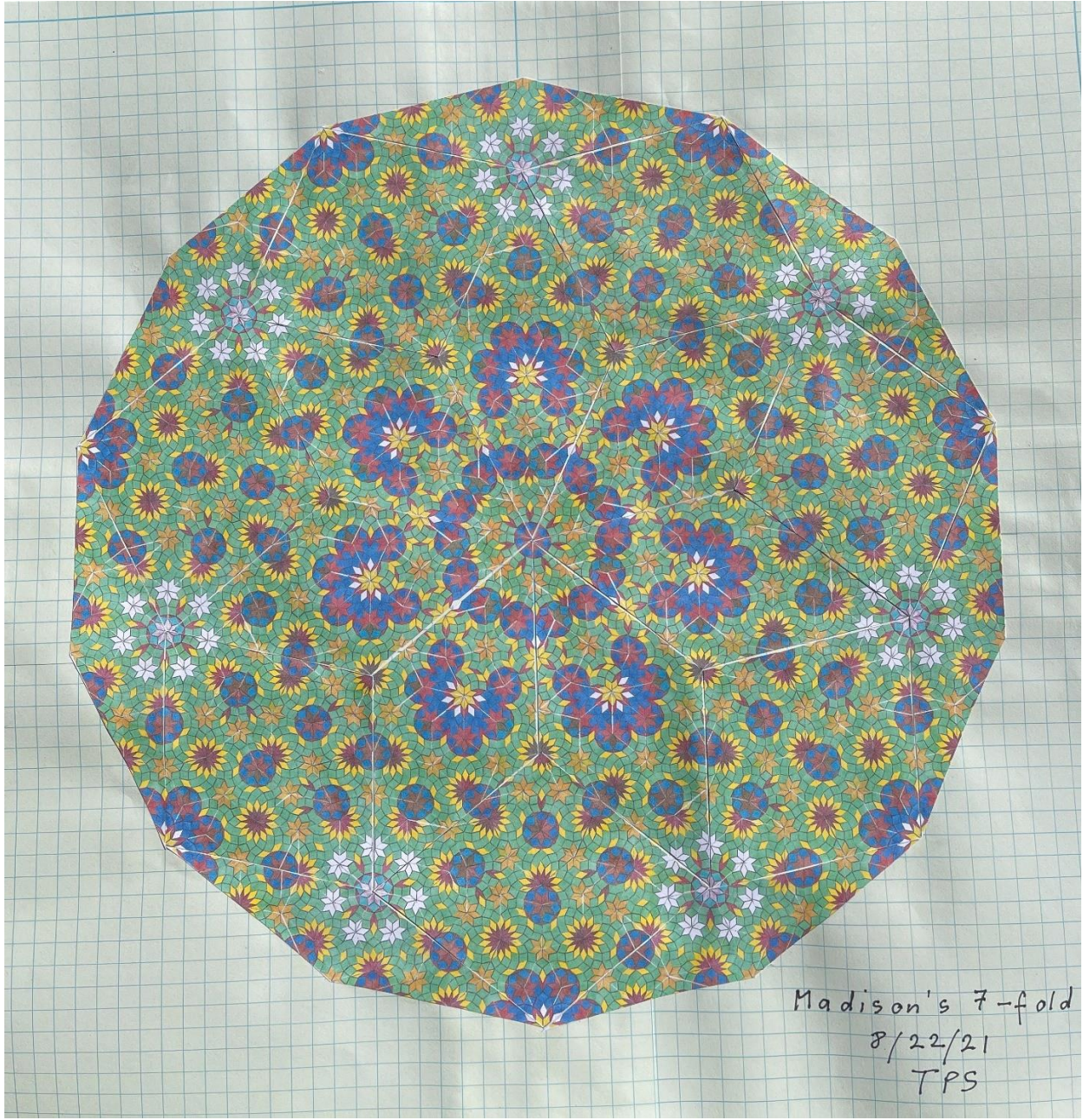

**Fig. 9:** A patch of Madison's 7-fold tiling. There are five differently colored 7-pointed stars and a 14-fold motif with 14 skinny rhombs, surrounded by yellow B tiles and green C tiles.

Alexey Madison also discovered another 7-fold tiling with six substitutions. It has an inflation factor of δ = 9.09783…

## Pattern #5, a 7-Fold Tiling with Seven Substitutions:

The Pattern #5 tiling has a substitution matrix

$$M^n = \begin{vmatrix} S_{4n}-S_{4n-1} & S_{4n} & S_{4n+1}-S_{4n} \\ S_{4n} & S_{4n+1}+S_{4n-1} & S_{4n+1} \\ S_{4n+1}-S_{4n} & S_{4n+1} & S_{4n}+S_{4n+1}-S_{4n-1} \end{vmatrix}$$

where $S_0 = 0$, $S_1 = 0$, $S_2 = 1$

$S_n = 2S_{n-1} + S_{n-2} - S_{n-3}$ ($n \geq 3$)

$S_{n=0,1,2...} = 0, 0, 1, 2, 5, 11, 25, 56, 126, 636, 1429, ...$

$S_n$ is a third-order Fibonacci series.

$$M = \begin{vmatrix} 3 & 5 & 6 \\ 5 & 9 & 11 \\ 6 & 11 & 14 \end{vmatrix}$$

is the substitution for the first generation. The inflation factor is $\delta = \rho^2$ (= 5.04892…), which is the square of the longer diagonal of a regular heptagon.

Figure 10 shows a new solution with seven substitutions. To be distinct from the Madison tiling, I started out with a different skinny red tile ($A_1$ tile). The inflated $A_1$ tile then has again red tiles at the acute vertices. It follows, that if a 14-fold star is inflated, it will always have a red 14-fold star at its center. Madison would label such a vertex as C, a constant vertex that does not change. The rosette in Figure 12 displays this characteristic. A second (deep purple) $A_2$ tile was needed to complete the inflation of $A_1$ in the first generation (the Madison tiling has only one A tile).

On the other hand, only three B tiles were needed. As a result, there are only three colored 7-fold stars (versus five in Madison's tiling). The rosettes in Figure 11 of these three stars show that the center of a white star turns into a yellow star in the next generation. Yellow turns to brown, and brown reverts to white again. Such a vertex is cyclical.

The background is dominated by the green $C_1$ tile. It is interesting to compare it with the dominant green $R_7$ tile in the Madison tiling. In both cases, the tiles near the vertices are partial stars, but with different base tiles.

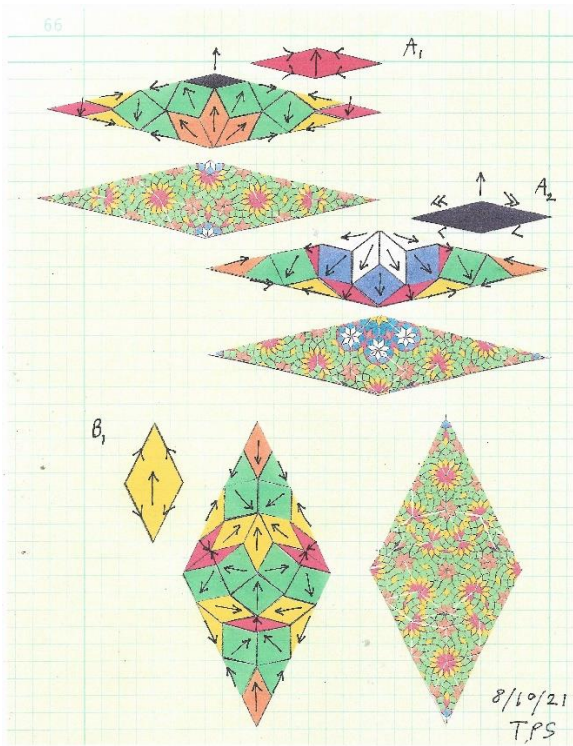
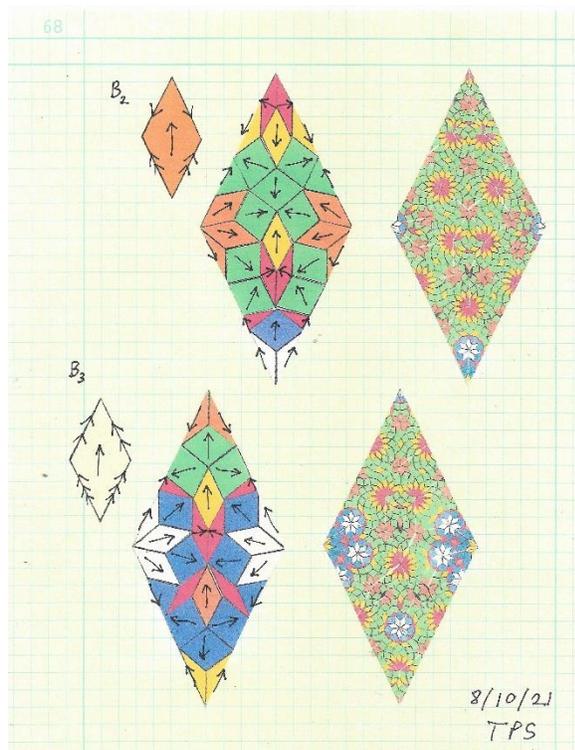
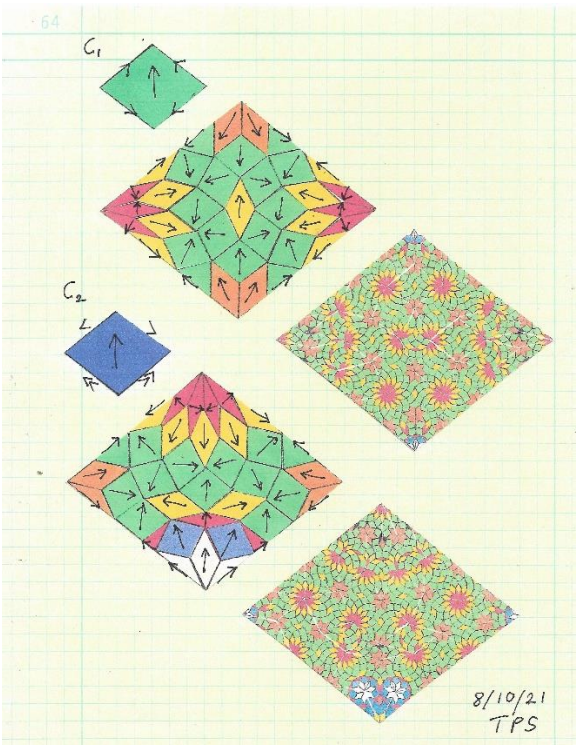

**Fig. 10:** Substitution and matching rules for the 7-fold Pattern #5 tiling. The background is dominated by the green $C_1$ tile. There are three distinct edges marked by arrows. The inflation sequences are:

$\delta = 1 + \Phi + \psi + 1$ (single arrow)

$\delta = \Phi + \gamma + 1 + \Phi$ (double arrow)

$\delta = \Phi + 1 + \psi + 1$ (triple arrow)

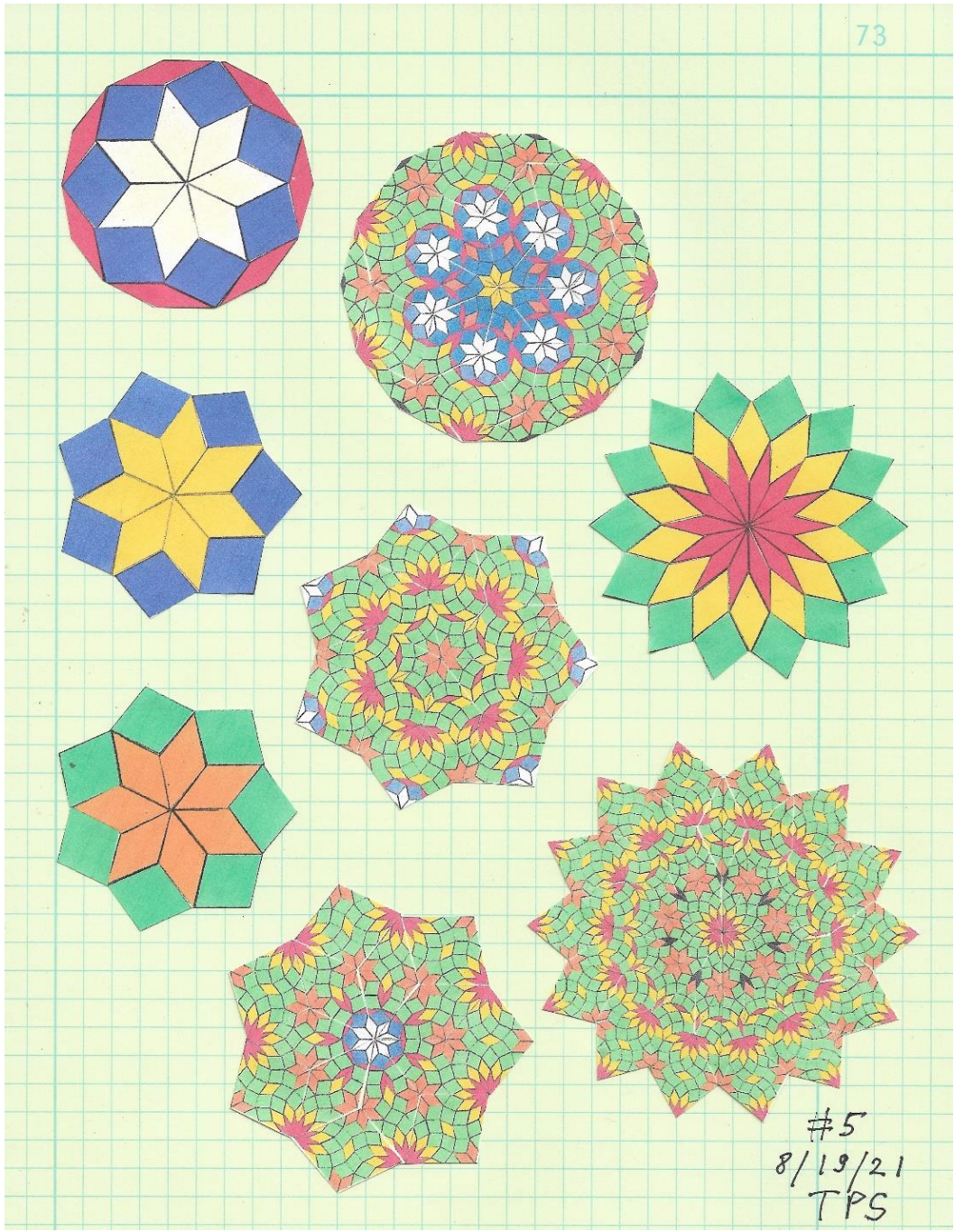

**Fig. 11:** The Pattern #5 has four main motifs that reoccur in a quasiperiodic manner. The first inflation is shown for each rosette. The 7-fold stars change colors at the center in a cyclical way. The 14-fold star repeats itself at the center and stays constant.

## Starry Night

Returning to the star theme throughout this paper, I tried to create a celestial rendition of Pattern #5 with a blue background and yellow 7-pointed stars. Conveniently, a purple $B_1$ tile also adds to the background while it also appears as a purple star, surrounded by the lighter $C_2$ tiles. I colored $C_2$ light green in connection with the purple star but left it white inside the tetradecagon motif. It should be considered the same base tile.

The patch shown in Figure 15 has all the main motifs mentioned before: three 7-fold stars in different colors, the red one inside a tetradecagon, and 14-fold stars. There are also partial 14-fold stars. The orange skinny tiles, together with purple tiles, also form ribbons. If the brilliant yellow stars are light matter, the ribbons could be dark matter.

The celestial theme was partially motivated by starry decorations that the author remembered from baroque ceilings and paintings. He is currently searching for some examples. Starry Night is, of course, the title of a famous painting by Vincent Van Gogh painted in 1889. Another iconic rendition of a starry sky is the Flammarian engraving of 1888.

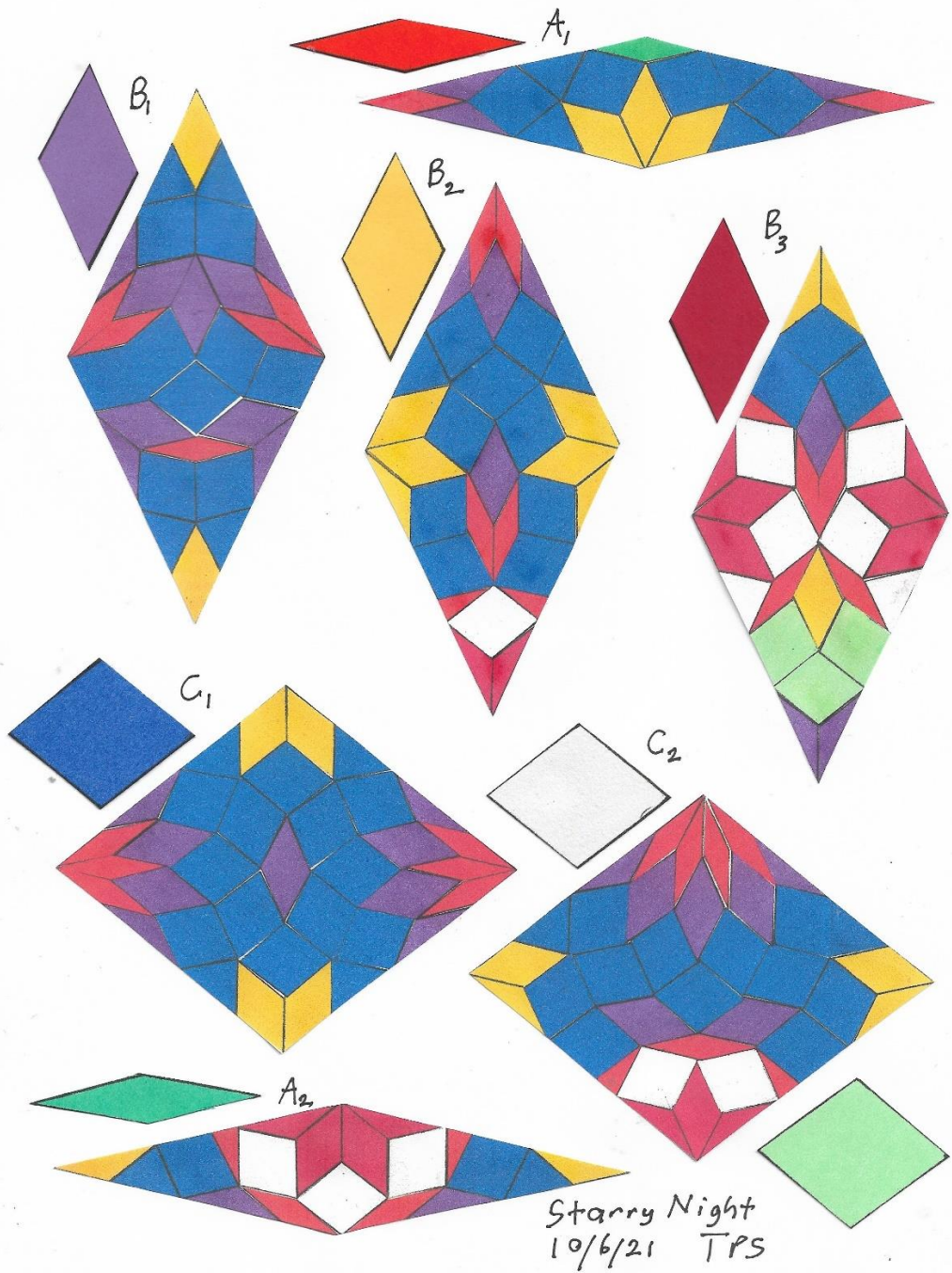

**Fig. 12:** Starry Night – The first generation of Pattern #5 in a different color scheme to create a celestial background in blue and purple and dominant yellow 7-pointed stars. The white and light green $C_2$ tiles follow the same substitution rule.

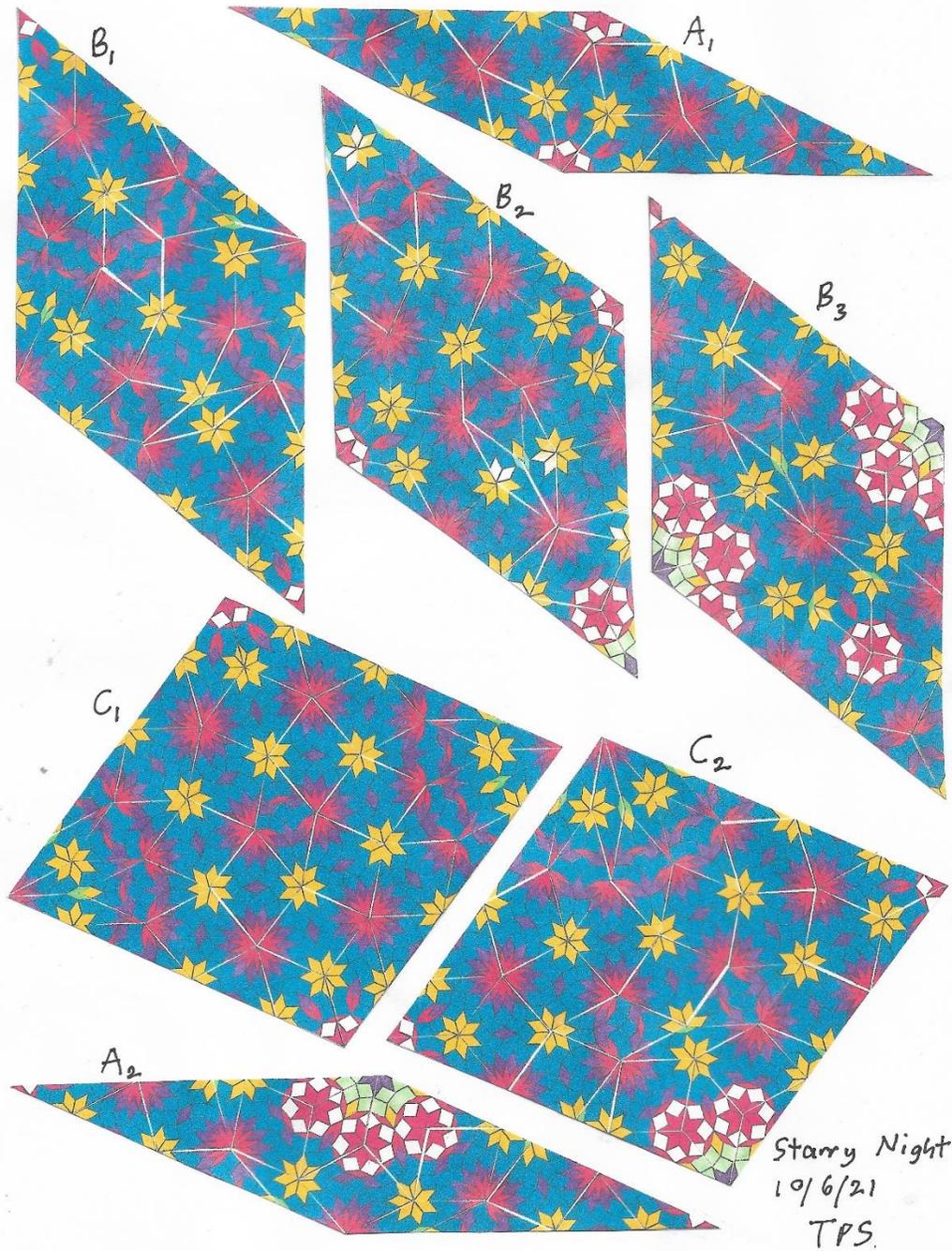

**Fig. 13:** Starry Night – The second generation of Pattern #5 with a celestial theme.

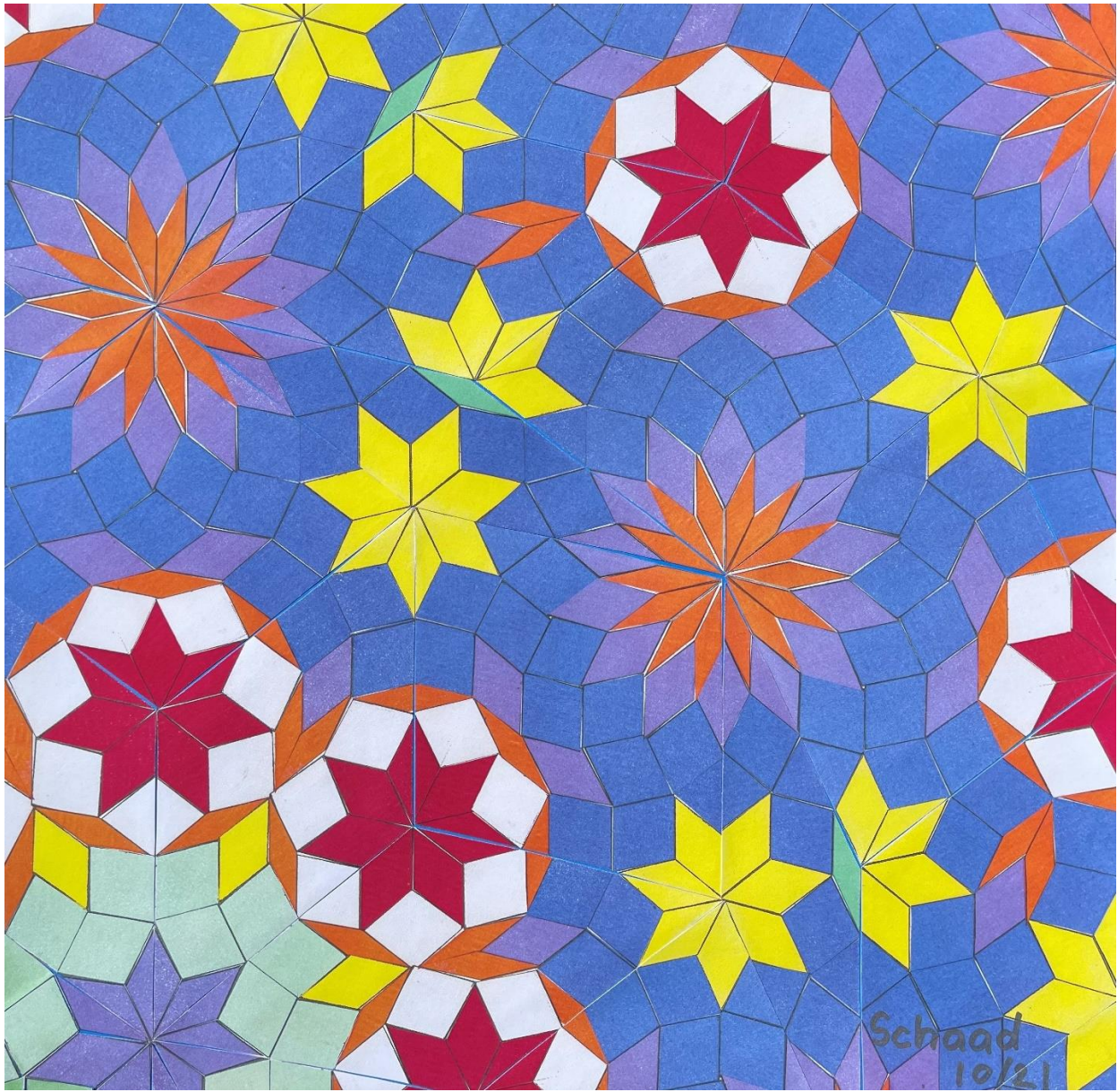

**Fig. 14:** Starry Night – a patch of the rhombic 7-fold tiling #5 with the background tiles $C_1$ and $B_1$ in blue and purple. There are four main motifs: three 7-fold stars in different colors and a 14-fold star of skinny rhombs in orange. In this view, there are also two odd pieces: an irregular hexagon with a purple tile and two skinny rhombs, and an unfinished yellow star with four points, capped with a green skinny rhomb. Whenever a green tile shows up in a larger skinny tile, there is an ambiguity if the larger tile should be flipped by 180 degrees.

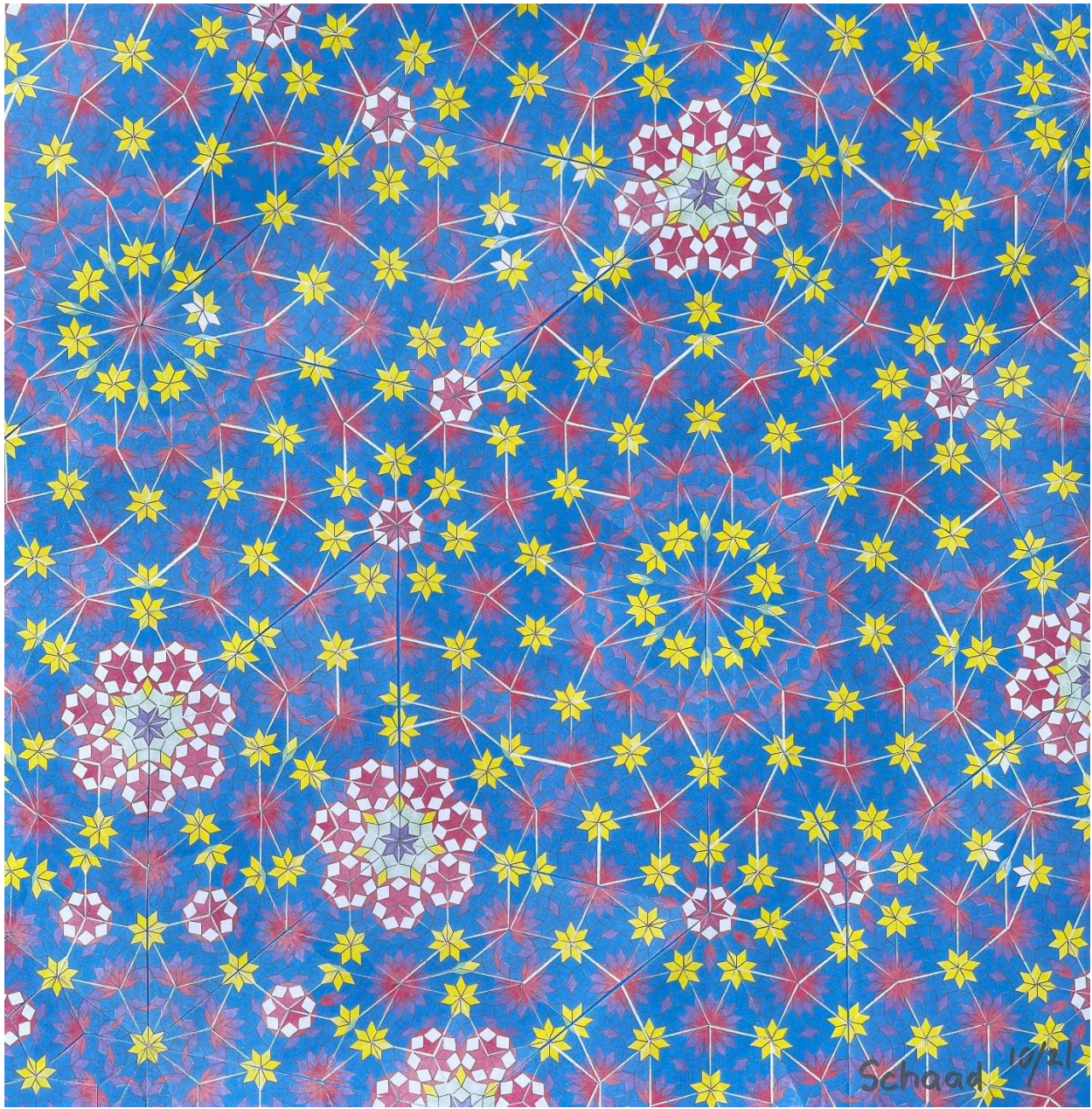

Fig. 15: Starry Night – The same patch of Pattern #5 as in the previous figure but inflated into the next generation. The three colored stars of Figure 15 have maintained their 7-fold rotational symmetry but have now morphed into a different color at its center.  The two 14-fold stars maintain another 14-fold star at its center, surrounded by a combination of 7-fold stars, partial stars, and irregular hexagons. Further out, those centers are surrounded by 14 yellow stars. The partial stars closer to the center break the local 7-fold rotational symmetry. Could there be matching rules that preserve it?